\documentclass[11pt, letterpaper]{article}

\usepackage{lscape,epsfig}
\usepackage{fancyhdr, color}
\usepackage{amsmath, amsfonts}
\usepackage{natbib}

\usepackage[english]{babel}

\usepackage{enumerate}
\newcommand{\cip}{\mbox{$\perp\!\!\!\perp$}}
\newtheorem{defn}{Definition}[section]
\newtheorem{lemma}[defn]{Lemma}
\newtheorem{theorem}[defn]{Theorem}

\newtheorem{assumption}[defn]{Assumption}
\newtheorem{example}[defn]{Example}

\title{Optimal estimation of coarse structural nested mean models with application to initiating ART in HIV infected patients}

\author{Judith J. Lok\\
Department of Mathematics and Statistics, Boston University\\ 111 Cummington Mall, Boston, Massachusetts 02215, U.S.A.\\ Email: jjlok@bu.edu}

\begin{document}

\maketitle

\section*{Abstract} 

Coarse structural nested mean models are used to estimate treatment effects from longitudinal observational data. Coarse structural nested mean models lead to a large class of estimators. It turns out that estimates and standard errors may differ considerably within this class. We prove that, under additional assumptions, there exists an explicit solution for the optimal estimator within the class of coarse structural nested mean models. Moreover, we show that even if the additional assumptions do not hold, this optimal estimator is doubly-robust: it is consistent and asymptotically normal not only if the model for treatment initiation is correct, but also if a certain outcome-regression model is correct.
We compare the optimal estimator to some naive choices within the class of coarse structural nested mean models in a simulation study. Furthermore, we apply the optimal and naive estimators to study how the CD4 count increase due to one year of antiretroviral treatment (ART) depends on the time between HIV infection and ART initiation in recently infected HIV infected patients. Both in the simulation study and in the application, the use of optimal estimators leads to substantial increases in precision.\\

\noindent {\bf Keywords:} Causal inference; Doubly robust estimation; HIV/AIDS; Longitudinal data; Observational studies; Structural nested models

\section{Introduction}

The effect of time-dependent treatments is often estimated from observational data, since clinical trials where treatment is repeatedly randomized are not common. Estimating treatment effects from observational data is more difficult than from clinical trials. Since treatment was not randomized, patients receiving treatment typically have different pre-treatment character-istics than patients not receiving treatment, leading to confounding by indica-tion. If patients
with a worse prognosis were treated more often, a naive analysis would lead to underestimation
of the treatment effect. It could even reverse the sign of the treatment effect, as illustrated for
time-dependent treatments in an HIV example in \cite{Victor}.

Several approaches exist to estimate treatment effects from longitudinal observational data.
Time-dependent coarse structural nested mean models (coarse SNMMs) describe the effect
of time-dependent treatments, conditional on patient characteristics at the time of treatment initiation. Coarse SNMMs model the mean difference between the outcome with treatment initiated at time
$m$, versus never, on the outcome measured at a later time $k$, given a patient was not treated until time
$m$, and given a patient’s covariate history at time $m$ (\citeauthor{Victor} \citeyear{Victor}). Previously,
\cite{Rnoncomp} introduced coarse SNMMs for outcomes measured at the
end of a study, for trials with noncompliance. Earlier SNMMs studied in \cite{Aids}, \cite{smoke}, \cite{SNart}, and \cite{comm} describe the effect of one treatment dosage conditional on patient characteristics just prior to this
dosage. Estimating the effect of a longer duration of treatment from these quantities requires modeling the distribution of covariates given past treatment and covariate history,
possibly leading to bias (\citeauthor{R86} \citeyear{R86}, \citeyear{R87b}, \citeyear{R89}, \citeyear{comm}), or even incompatibility of the
different assumptions (\citeauthor{comm} \citeyear{comm}). Coarse structural nested mean models are thus useful to estimate
the effect of a treatment that is initiated once and then never stopped, since they directly
estimate the effect of multiple treatment dosages.

Marginal structural models (\citeauthor{MSM1} \citeyear{MSM1}, \citeauthor{MSM2} \citeyear{MSM2}) are another class of models to estimate the effect of time-dependent treatments from observational data. Marginal structural models estimate how treatment effects depend
on baseline, but not time-dependent, covariates. \cite{R99} provides a detailed comparison of marginal structural models and structural nested models. Q-learning (\citeauthor{ChakrabortyMurphy2014} \citeyear{ChakrabortyMurphy2014}, \citeauthor{NahumMurphy2012} \citeyear{NahumMurphy2012}) also estimates treatment effects from time-dependent observational data (and SMART trials, where treatments are repeatedly randomized). Q-learning builds on classical dynamic programming: 1. estimate the optimal last treatment strategy; 2., with backwards
induction, the optimal treatment strategy at each time point given the optimal treatment
strategy will be followed onwards. Q-learning is thus aimed at estimating optimal
treatment regimes, and does not focus on the effect of a fixed duration of treatment, or of treatment
initiation at different time points. As
traditional SNMMs, \cite{Almi2013} estimate the effect of a last blip of
treatment conditional on a covariate and treatment history, and how
this effect depends on potential effect modifiers, while allowing for confounders
one might not want to estimate the effect modification of. They do not cover optimal estimation.

SMART trials (\citeauthor{ChakrabortyMurphy2014} \citeyear{ChakrabortyMurphy2014}, \citeauthor{Almi2014} \citeyear{Almi2014}) determine the best time-dependent treatment strategies using randomized trials. In SMART trials, baseline and subsequent treatment decisions are randomized, with randomizaton options determined by prior outcomes. \cite{LiuKosorok} combine outcome weighted learning with Q-learning to solve the optimal
treatment strategy from SMART trials. They use support vector machines after a reformulation
of the expected utility under any treatment strategy that leads to a convex optimization problem.

Various methods have been proposed to efficiently estimate the effect of point treatments. \cite{Newey1993} describes efficient estimation under conditional moment restrictions. Our identifying
assumption leads to a conditional mean being independent of a time-dependent treatment
variable, which cannot be written in terms of such a conditional moment restriction, so the theory
from \cite{Newey1993} cannot be applied to derive efficient estimators of coarse structural nested
models. There are similarities between our work and \cite{Newey1993} in that we are deriving an
optimal “instrumental variable” in the sense of the econometrics literature: a function of the
covariates which leads to the smallest asymptotic variance. Targeted maximum likelihood estimation
(TMLE: \citeauthor{LaanRubin} \citeyear{LaanRubin}) is a method for doubly robust and efficient estimation
of parameters in several semiparametric and nonparametric settings. Later targeted maximum
likelihood estimation theory (\citeauthor{LaanGruber} \citeyear{LaanGruber}, \citeauthor{Schnitzer} \citeyear{Schnitzer}, \citeauthor{Petersen} \citeyear{Petersen}) builds on theory from \cite{Bang}. Targeted maximum likelihood estimation
has not been developed for structural nested models. \cite{Tsiatis} covers semiparametric efficient
estimation, but not of structural nested models. \cite{Tsiatis} focuses on missing data
problems, which could potentially be adapted easier to marginal structural models (\citeauthor{MSM1} \citeyear{MSM1}, \citeauthor{MSM2} \citeyear{MSM2}), which use inverse probability of treatment weighting. Inverse
probability of treatment weighting methods are also advocated in \cite{ChakrabortyMurphy2014}. Matching (\citeauthor{Guido2004} \citeyear{Guido2004}, \citeauthor{Rosenbaumbook} \citeyear{Rosenbaumbook}) has been restricted to
baseline treatments, not time-dependent treatments, and the most commonly used version with
a fixed number of matches per individual is not efficient (\citeauthor{AbadieImbens2006} \citeyear{AbadieImbens2006}).  \cite{GutmanRubin} describe the effect of a point treatment on a binary outcome $Y$ in the presence
of one continuous confounder $X$. They use $EY = EE[Y|X]$ for $Y$ the outcome under no treatment and for $Y$ the
outcome under treatment (see also \citeauthor{Bang} \citeyear{Bang} Section~2.1). \cite{GutmanRubin} estimate the conditional
expectations separately for treatment and control with splines, and use multiple imputation
instead of the estimated conditional means. This method is not doubly robust or efficient (\citeauthor{Bang} \citeyear{Bang}). \cite{Hahn} proposed efficient methods to estimate the effect of
point treatments from observational data. His proposed imputation methods are hard to generalize
to time-dependent treatments, where if there are $K$ potential treatment times, the number
of missing potential outcomes per patient is $2^K- 1$. In addition, it seems easier to specify the
treatment initiation models than the models for the conditional expectations of the outcomes.  The theory of \cite{Chernozhukov2015} Section~3 on optimal variance does not apply to coarse structural nested mean models, since, as we will show, coarse structural nested mean models lead to a continuum of (orthogonal) unbiased estimating equations for the causal parameter of interest. 

In contrast to most of the literature on efficient estimation of treatment effects from observational
data (\citeauthor{Guido2004} \citeyear{Guido2004} provides an overview), this article focuses on optimal estimation of
the effect of a time-dependent treatment. In addition, coarse structural nested mean models can estimate how
the treatment effect depends on time-dependent pre-treatment characteristics, which is not the
focus of most of the literature on efficient estimation of treatment effects.

We provide a description of the estimation methods and the assumptions needed to
consistently estimate coarse structural nested mean models (coarse SNMMs) with an outcome measured over time. Detailed proofs are in \cite{Victor}, \cite{comm}, and \cite{RRS}. The main focuses of this article however
are double robustness and optimal estimation.

\cite{Rnoncomp} and, for time-dependent outcomes, \cite{Victor} derived a large
class of estimating equations for coarse SNMMs, all leading to consistent,
asymptotically normal estimators for the treatment effect. It turns out that both estimates
and standard errors may depend considerably on the choice of estimating equations within this
large class. This motivates the current article, which derives, under extra conditions, an optimal choice of estimators within a class that includes the estimators from \cite{Victor}. This optimal estimator leads to the smallest possible asymptotic variance. It is also doubly robust: it is
consistent and asymptotically normal not only if the model for treatment initiation is correctly
specified, but also if a certain outcome-regression model is correctly specified. The optimal estimator combines weighting and regression, and is therefore a mixed method (compare
with mixed methods for point exposures described in \citeauthor{Guido2004} \citeyear{Guido2004}).

Also without the extra conditions needed for optimality and in the presence of censoring,
this optimal estimator is doubly robust; it just may not be optimal in such settings. As shown in Application
Section~\ref{results}, the optimal estimator may perform considerably better than arbitrarily choosing an estimator within
the class of coarse SNMMs.

We implemented the proposed optimal estimator and compare it to other estimators for coarse
SNMMs, both in a simulation study and in an HIV application.We estimate how
the effect of one year of ART, the current standard of care for HIV infected patients, depends
on the time between the estimated date of HIV infection and ART initiation. The application
includes correction for informative censoring, and bootstrap confidence intervals. Consistency
of the bootstrap for all our estimators is proven in Web-appendix~A under regularity conditions.

Optimality for coarse SNMMs is simpler than optimality for more traditional
SNMMs (see \citeauthor{comm} \citeyear{comm}, \citeauthor{Laan} \citeyear{Laan}),
because coarse SNMMs avoid the need for accumulating the effects over
time. This article therefore
includes an accessible illustration of the steps involved in calculating optimal estimators.
Our Web-Appendix explicitly proves all our results for time-dependent coarse SNMMs, thus providing a self-contained example.


\section{Setting and notation}
\label{setting}

Initially, all patients are assumed to be followed at the same times $0,1,\ldots,K+1$, with $0$ indicating the baseline visit. The assumption that there is no censoring
due to loss-to-follow-up is relaxed in Section~\ref{steps}. $Y_{k}$ is the outcome at time $k$. $\overline{Y}_{k}=(Y_0,Y_1,\ldots,Y_k)$ is the outcome history until time $k$. $A_k$ is the treatment at time $k$. We investigate the effect of a binary treatment $A_k$, which is either given ($A_k=1$) or not
($A_k=0$) at each time $k$. $\overline{A}_{k}$ is the treatment
history until time $k$. We only consider the impact of initiating
treatment, and do not consider issues of treatment interruption or
compliance: our analysis follows the intention-to-treat
principle common in randomized trials. Thus, $A_k=0$ until treatment is initiated, and $A_k=1$ thereafter. $T$ is the time treatment was actually initiated, with $T=K+1$ if treatment was never initiated. Similarly, $L_k$
are the covariates at time $k$, including the outcome $Y_k$ at time $k$, and $\overline{L}_{k}$ is the
covariate history until time $k$. $\overline{{\cal L}}_{k}$ is the
space in which $\overline{L}_{k}$ takes its values. We assume that at
each visit time $k$, treatment decisions $A_k$ are made after $L_k$
is measured and known. We suppress patient-level notation (such as the subscript $i$ often used to indicate individual $i$).

$Y^{(\emptyset)}_k$ and $\overline{Y}^{(\emptyset)}_k$ are the counterfactual, not always measured, outcome and outcome history at time $k$ under the treatment regime ``no treatment''.  $L^{(\emptyset)}_k$ and $\overline{L}^{(\emptyset)}_k$ are the covariates and covariate history at time $k$ under the treatment regime ``no treatment''. $Y_k^{(m)}$ is the outcome at time $k$ under the treatment regime ``start treatment at time $m$''. We assume that observations and counterfactuals of the different patients are independent and identically distributed (\citeauthor{Rubin} \citeyear{Rubin}).

\section{Time-dependent coarse structural nested mean models}\label{smodel}

Our model for treatment effect is similar to that in \cite{Rnoncomp},
but differs in that we allow a time-dependent outcome, as in
\cite{Victor}:
\begin{defn}\label{gamma} \emph{(TIME-DEPENDENT COARSE SNMMS).}\\ For $k=m,\ldots, K+1$,
\begin{equation*}
\hspace*{3cm}\gamma_k^{m}\left(\overline{l}_m\right)=E\left[Y_{k}^{(m)}-Y_{k}^{(\emptyset)}\mid \overline{L}^{(\emptyset)}_{m}=\overline{l}_m,T=m\right].
\end{equation*}
\end{defn}
$\gamma_k^{m}\left(\overline{l}_m\right)$ is the expected difference, for patients whose treatment started
at time $m$ with covariates $\overline{l}_m$, of 
the outcome at time $k$ had the patient started treatment at time
$m$, and the outcome at time $k$ had the
patient never started treatment. It is the effect of treatment between times $m$ and $k$.

We assume a parametric model for the treatment effect $\gamma$, leading to a semiparametric setting. $\gamma$ models the expected difference between two counterfactual outcomes. How any of these outcomes depends on a patient's covariates at time $m$ is not specified. For example, it could be that, for $k>m$ and with $g_k(\overline{L}_m)$ any function of $\overline{L}_m$,
\begin{equation*}
Y_k^{(\emptyset)}=g_k(\overline{L}_m)+\epsilon^{(\emptyset)}_{k}, \hspace{1cm} Y_k^{(m)}=g_k(\overline{L}_m)+\left(\psi_1+\psi_2m+\psi_3m^2\right)(k-m)+\epsilon^{(m)}_{k}
\end{equation*}
(note that $k-m$ is the treatment duration until time $k$), with\\ $E[\epsilon^{(m)}_{k}\mid\overline{L}_m,T=m]=0$ and $E[\epsilon^{(\emptyset)}_{k}\mid\overline{L}_m,T=m]=0$. In this example, $g_k(\overline{L}_m)$ is left completely unspecified, and
\begin{equation*}
\gamma^{m}_{k,\psi}\left(\overline{l}_m\right)=
\left(\psi_1+\psi_2m+\psi_3m^2\right)(k-m)1_{\left\{k>m\right\}}.
\end{equation*}
So-called rank preservation holds if $Y_k^{(m)}-Y_k^{(\emptyset)}=\gamma_k^{m}\left(\overline{L}_m\right)$, but rank preservation has been argued to not hold in many practical applications (\citeauthor{Enc} \citeyear{Enc}, \citeauthor{MimAn} \citeyear{MimAn}). A coarse SNMM is a model for the effect of the treatment on the treated, as discussed in e.g.~\cite{Guido2004}. If $\overline{L}$ contains enough information to make those treated and those untreated at time $m$ comparable given $\overline{L}_m$, a structural nested mean model is also a model for the effect of the treatment on all patients with given covariates $\overline{L}_m$. The lack of distinction between the effect of the treatment on the treated and the effect of the treatment is due to the fact that structural nested models condition on pre-treatment covariates.

\begin{assumption} \label{model} \emph{(Parameterization of coarse structural nested mean model).} $\gamma^{m}_{k,\psi}\left(\overline{l}_m\right)$ is a correctly specified model for $\gamma_k^{m}\left(\overline{l}_m\right)$, with $\psi$ a parameter in ${\mathbb R}^p$ and $\psi_{*}$ the true parameter.
\end{assumption}
The following example (\citeauthor{Victor} \citeyear{Victor}) motivated this work on optimal estimation:
\begin{example}\label{exART} (Effect of antiretroviral treatment (ART) depending on the time between estimated date of HIV infection and treatment initiation in HIV infected patients).
We initially assume that
\begin{equation}\label{modelex}
\gamma^{m}_{k,\psi}\left(\overline{l}_m\right)=
\left(\psi_1+\psi_2m+\psi_3m^2\right)(k-m)1_{\left\{k>m\right\}},
\end{equation} with $(k-m)$ the treatment duration from
month $m$ to month $k$, is a correctly specified parametric model. Possibly, the mean treatment effect also
depends non-linearly on the treatment duration. In that case, one
could add non-linear terms such as $\psi_4(k-m)^21_{\left\{k>m\right\}}$, or additionally include a non-linear term
depending on time since infection, such as $(\psi_4+\psi_5m)(k-m)^21_{\left\{k>m\right\}}$.

The treatment effect $\gamma$ may also depend on pre-treatment covariates, such as log$_{10}$ HIV
viral load in the blood ($lvl$), resistance mutations, and the CD4 count. To incorporate pre-treatment covariates such as log$_{10}$ HIV viral load, one can extend the model by
including terms such as $\psi_4 lvl_m (k-m)1_{\left\{k>m\right\}}$. We chose some
function of $(k-m)$ because the treatment duration may be
predictive of its effect. If, for example, the treatment effect
depends only on the viral load at treatment initiation for the
first month of treatment, one might add terms such as $\psi_4 lvl_m
1_{\left\{k>m\right\}}$.
\end{example}

Following \cite{Aids} and \cite{comm}, we use the
propensity score (\citeauthor{prop} \citeyear{prop}), the prediction of treatment given the past, $p(m)= pr\left(A_m=1\mid \overline{A}_{m-1}=\overline{0},\overline{L}_{m}\right)
$, to
estimate the treatment effect $\gamma$. Henceforth, we assume that $p_\theta(m)$ is a
correctly specified model for $p(m)$, with $p(m)=p_{\theta_{*}}(m)$.
We typically estimate $\theta_*$ by maximum
partial likelihood.

\section{No unmeasured confounding and consistency}
\label{nuc}

As in \cite{Aids}, \cite{Enc}, \cite{MSM1}, \cite{SNart},
\cite{Rnoncomp}, and \cite{Victor}, to distinguish between treatment
effect and confounding by indication, we require the assumption of no unmeasured confounding. It states that information is
available on all factors that both: (1) influence treatment decisions
and (2) possibly predict a patient's prognosis with respect to the
outcome of interest. $Y^{(\emptyset)}_{k}$, the outcome at
time $k$ without treatment, reflects a patient's
prognosis with respect to the outcome of interest. If there is no unmeasured
confounding, treatment decisions at time $m$ ($A_m$) are
independent of this (not always measured) prognosis $Y^{(\emptyset)}_{k}$ ($k>m$), given
past treatment and covariate history $\overline{A}_{m-1}$ and
$\overline{L}_{m}$:
\begin{assumption}\label{intconf}\emph{(No unmeasured confounding -
    formalization)}. $A_m\cip Y^{(\emptyset)}_{k}
\mid \overline{L}_{m},\overline{A}_{m-1}$
for $k>m$, where $\cip$ means: is independent of
(\cite{Dawid}).
\end{assumption}

If
a patient is not treated until time $k$, there is no difference in
treatment between $Y_{k}$ and $Y^{(\emptyset)}_{k}$. In this and similar
cases, it is reasonable to assume that until time $k$, the observed outcomes 
and the outcomes without treatment would have been the same:
\begin{assumption} \label{cons}\emph{(Consistency).}
If $T\geq k$, $Y_{k}=Y^{(\emptyset)}_{k}$ and $\overline{L}_k=\overline{L}_k^{(\emptyset)}$. $Y^{(T)}=Y$.
\end{assumption}

\section{Estimation: unbiased estimating equations}
\label{est}

This section describes the estimation methods from \cite{Rnoncomp} and
\cite{Victor}. Proofs can be found in \cite{Victor}.
\begin{defn} \label{Hdef}
On $k>T$, define $H(k)=Y_{k}-\gamma_{k}^{T}\left(\overline{L}_{T}\right)$. On $k\leq T$, define $H(k)=Y_k$.
\end{defn}
\begin{example}\label{exART2} (Effect of antiretroviral treatment (ART) depending on the time between estimated date of HIV infection and treatment initiation in HIV infected patients).
In the setting of Example~\ref{exART}, on $k>T$,
$
H(k)=Y_k-\left(\psi_1+\psi_2 T +\psi_3 T^2\right)(k-T)
$.
\end{example}
For the true $\psi_{*}$, $H_{\psi_{*}}(k)=Y_k-\gamma_{k,\psi_{*}}^{T}\left(\overline{L}_T\right)=H(k)$. This so-called blipping off of the treatment effect generates a random variable $H(k)$ that mimics a counterfactual outcome:
\begin{theorem} \label{mim}\emph{(MIMICKING COUNTERFACTUAL OUTCOMES).} Under consistency assumption~\ref{cons}, for $m\leq K$ and $k\geq m$,
\begin{equation*}
E\left[H(k)\mid \overline{L}_m,\overline{A}_{m-1}=\overline{0},A_m\right]
=E\left[Y_{k}^{(\emptyset)}\mid \overline{L}_m,\overline{A}_{m-1}=\overline{0},A_m\right].
\end{equation*}
\end{theorem}

The idea behind estimation, similar to for example \cite{Aids} or
\cite{SNart}, is that under assumption of no unmeasured confounding~\ref{intconf}, given past treatment and covariate history, $Y^{(\emptyset)}_{k}$ does not help to
predict treatment changes. In other words, in the model for treatment changes, that is, the propensity score (\citeauthor{prop} \citeyear{prop}), $Y_k^{(\emptyset)}$ does not contribute.
\cite{Victor} proved that because of Theorem~\ref{mim},
the same
holds for $H_{\psi_{*}}\left(k\right)$, and that similar to \cite{ASarx,plss},
this leads to the
following theorem:
\begin{theorem}\emph{(UNBIASED ESTIMATING EQUATIONS).}\label{see}
  Suppose that consistency assumption~\ref{cons} and assumption of
  no unmeasured confounding~\ref{intconf} hold. Consider any
  $\vec{q}_{m}^{\;k}:\overline{{\cal L}}_{m}\rightarrow
  \mathbb{R}^{p}$, $m=0,\ldots,K$, $k>m$, which are
  measurable, bounded, and vector-valued.  Then
\begin{equation*}
E \left(\sum_{m=0}^{K}\sum_{k>m}
\vec{q}^{\;k}_{m}\left(\overline{L}_{m}\right) H\left(k\right) 1_{\overline{A}_{m-1}=\overline{0}}
\left\{A_m-p\left(m\right)\right\}\right)=0.
\end{equation*}
If furthermore $\gamma_\psi$ is correctly
specified (Assumption~\ref{model}) and $p_\theta\left(m\right)$
correctly specifies $p\left(m\right)$, then
\begin{equation*}
P_n\left(\sum_{m=0}^{K}\sum_{k>m} \vec{q}^{\;k}_{m}\left(\overline{L}_{m}\right) H_{\psi}\left(k\right) 1_{\overline{A}_{m-1}=\overline{0}}
\left\{A_m-p_\theta\left(m\right)\right\}\right)=0,
\end{equation*}
stacked with
the estimating equations for $\theta_{*}$, with $P_n$ the empirical
measure $P_n(X)=1/n \sum_{i=1}^n X_i$, are unbiased estimation
equations for $(\psi,\theta)$. The $\vec{q}^{\;k}_{m}$ here are allowed to depend on $(\psi,\theta)$, as long as they are measurable and bounded for
$\left(\psi_{*},\theta_{*}\right)$.
\end{theorem}

\noindent For identifiability of the estimator, one needs as many
estimating equations as parameters, in this case by choosing the
dimension of $\vec{q}$. Including $k\leq m$
does not help in these estimating equations, since for those $k$, on
$\overline{A}_{m-1}=\overline{0}$, $H_\psi(k)=Y^{(\emptyset)}_{k}=Y_k$,
which is part of $\overline{L}_{m}$ and therefore generates a term
with expectation $0$ regardless of $\psi$.

If $\gamma$ is linear in $\psi$, this approach leads to
a linear restriction on $\psi$ once the parameter $\theta$ has been
estimated, and thus to a closed form expression for $\hat{\psi}$.

\section{Optimal estimating equations}\label{opt}

\subsection{Assumptions and restrictions on the estimating equations}

The vast literature on unbiased estimating equations, see
e.g.~\cite{Vaart} Chapter~5, indicates that under regularity
conditions, $\hat{\psi}$ is consistent and asymptotically normal for
any choice of $\vec{q}$ with the same dimension as $\psi$. Under those
regularity conditions, the asymptotic variance of an estimator
$\hat{\psi}$ that is a zero of $P_nG_\psi$, with $EG_{\psi_{*}}=0$, is
equal to
\begin{equation}\left\{E\left(\left.\partial/\partial \psi\right|_{\psi_{*}}G_\psi\right)\right\}^{-1}  E\left(G_{\psi_{*}}G^{\top}_{\psi_{*}}\right)  \left\{E\left(\left.\partial/\partial \psi\right|_{\psi_{*}}G_\psi\right)\right\}^{-1 \top},\label{asvar}
\end{equation}
because
\begin{equation}
n^{1/2}\bigl(\hat{\psi}-\psi_{*}\bigr)=\left\{E\left(\left.\partial/\partial \psi\right|_{\psi_{*}}G_\psi\right)\right\}^{-1} n^{1/2}P_n\bigl(G_{\psi_{*}}\bigr)+o_P(1).\label{expansion}
\end{equation}

In the following we restrict to estimators that satisfy the regularity
conditions needed for (\ref{expansion}) and thus (\ref{asvar}) to
hold. Among such estimators, this
article derives the optimal choice of $\vec{q}$. The remainder of this article also assumes assumption of no unmeasured confounding~\ref{intconf}, consistency
assumption~\ref{cons}, and correct specification of coarse structural nested mean model
assumption~\ref{model}. Web-Appendix~A provides proofs of all theorems and lemmas. Web-Appendix~A shows, under regularity conditions, consistency of the bootstrap for all estimators considered.

\subsection{Doubly robust estimators lead to increased precision}
\label{DRsec}

\begin{defn} \label{Gstar} Let 
\begin{equation*}
G(\psi,\theta,q)=
\sum_{m=0}^{K}\sum_{k>m} \vec{q}^{\;k}_m\left(\overline{L}_{m}\right)H_{\psi}(k)1_{\overline{A}_{m-1}=\overline{0}}\left\{A_m-p_{\theta}(m)\right\},
\end{equation*}
\begin{eqnarray*}
\lefteqn{G^*\left(\psi,\theta,q\right)=\sum_{m=0}^{K}\sum_{k>m} \vec{q}^{\;k}_m\left(\overline{L}_{m}\right)\left\{H_{\psi}(k)-E\left[H_{\psi}(k)\mid \overline{L}_m,\overline{A}_{m-1}=\overline{0}\right]\right\}}\\
&&\;\hspace*{7cm}
1_{\overline{A}_{m-1}=\overline{0}}\left\{A_m-p_{\theta}(m)\right\}.
\end{eqnarray*}
\end{defn}

First, we consider estimating equations for $\psi$ when $\theta_{*}$ and thus the propensity score is known:
\begin{theorem} \label{1or3} \emph{(REPLACEMENT OF ESTIMATING EQUATIONS BY\\ MORE EFFICIENT ONES).} Under no unmeasured confounding assumption~\ref{intconf}, consistency assumption~\ref{cons} and the usual regularity conditions for the sandwich estimator for the variance, based on equation~(\ref{asvar}), to hold,\\ $P_n\bigl(G^*\left(\psi,\theta_{*},q\right)\bigr)=0$
are unbiased estimating equations which, for given $\vec{q}$, lead to a smaller asymptotic variance of $\hat{\psi}$ than the estimating equations
$P_n\left(G(\psi,\theta_{*},q)\right)=0$.
\end{theorem}
The equations with $G^*$ are not
true estimating equations, because their specification depends on the parameter $\psi$ of interest, through
the conditional expectation of $H_\psi$. We will return to this issue
later. Theorem~\ref{mim} facilitates specifying the model for the conditional
expectation of $H$.

In practice, $\theta_{*}$ will usually be unknown and has to be
estimated. For the more efficient estimators based on $G^*$ of
Theorem~\ref{1or3}, estimating $\theta_{*}$ does not change the
asymptotic variance of $\hat{\psi}$. This result is similar to
Proposition~1 in \cite{RR95} for a missing data
problem:
\begin{theorem} \label{thetahat} Replacement of $\theta_{*}$ by $\hat{\theta}$ from a correctly specified pooled logistic regression model fitted by maximum partial likelihood, leading to $\hat{\psi}$ that solves $P_n\bigl(G^*\bigl(\psi,\hat{\theta},q\bigr)\bigr)=0$, leads to the same asymptotic variance for $\hat{\psi}$ as the estimator for $\psi_*$ that solves $P_n\bigl(G^*\left(\psi,\theta_{*},q\right)\bigr)=0$.
\end{theorem}
For the estimators solving $P_n\left(G(\psi,\theta_{*},q)\right)=0$, estimating $\theta_{*}$ may change the asymptotic variance. It usually reduces the asymptotic variance, as was also seen in \cite{optopt} and \cite{ASarx} for different  structural nested models. However, the resulting estimator is never more efficient than its doubly robust counterpart:
\begin{theorem} \label{3or1e} For $q$ fixed, replacement of $\theta_{*}$ by $\hat{\theta}$ from a correctly specified pooled logistic regression model fitted by maximum partial likelihood does not make the estimator for $\psi_{*}$ which solves $P_n\bigl(G\bigl(\psi,\hat{\theta},q\bigr)\bigr)=0$ more efficient than the estimator which solves $P_n\bigl(G^*\bigl(\psi,\hat{\theta},q\bigr)\bigr)=0$.
\end{theorem}
As for SNMMs (\citeauthor{optopt} \citeyear{optopt}, 3.10 page~23),
estimators resulting from $G^*$ are also doubly robust:
\begin{theorem} \label{DR} \emph{(DOUBLE ROBUSTNESS).} The estimator $\hat{\psi}$ which solves
$P_n\bigl(G^*\bigl(\psi,\hat{\theta},q\bigr)\bigr)=0$ is doubly robust: stacked with the estimating equations for $\theta$, these estimating equations are unbiased for $\psi$ if either $p_\theta$ or $E\left[H_{\psi}(k)\mid \overline{L}_m,\overline{A}_{m-1}=\overline{0}\right]$ is correctly specified. Thus, $\hat{\psi}$ is consistent and asymptotically normal if either of these models is correctly specified.
\end{theorem}
It follows that both for robustness and for efficiency, the estimators with $G^*$ are preferred.

\subsection{A theorem that guarantees optimality of estimators}
\label{optpropsec}

The following theorem, a consequence of Theorem~5.3 from \cite{Newey}, gives a sufficient criterion under which $\vec{q}^{\;opt}$ is optimal within our two classes of estimating equations, described by $G$ and $G^*$:
\begin{theorem}\label{optprop}\emph{(SUFFICIENT OPTIMALITY CRITERION WITH $\theta_{*}$ KNOWN).}
For both $G$ and $G^*$: if $\vec{q}^{\;opt}$ satisfies
\begin{equation}\label{opteqn}
E\left(\left.\partial/\partial \psi\right|_{\psi_{*}}G\left(\psi,\theta_{*},q\right)\right)
=E\left(G\left(\psi_{*},\theta_{*},q\right)G\left(\psi_{*},\theta_{*},\vec{q}^{\;opt}\right)^{\top}\right),
\end{equation}
then no other $\vec{q}$ satisfying our regularity conditions within this class leads to an estimator for $\psi_{*}$ with a smaller asymptotic variance than $\vec{q}^{\;opt}$. The estimator resulting from $\vec{q}^{\;opt}$ has asymptotic variance equal to the inverse of $E\bigl(G\left(\psi_{*},\theta_{*},\vec{q}^{\;opt}\right)G\left(\psi_{*},\theta_{*},\vec{q}^{\;opt}\right)^\top\bigr)$. There is a unique (in $L_2(P)$-sense) optimal solution to equation~(\ref{opteqn}) within this class of estimating equations.
\end{theorem}

\subsection{Explicit expression for optimal estimating equations for coarse SNMMs with a time-varying outcome}
\label{finaloptsec}

This section finds the optimal $\vec{q}$ under the following condition:
\begin{assumption}\label{HH} \emph{(Homoscedasticity).} For $0\leq m\leq K$ and $k,s>m$,
${\rm cov}\left[H(k),H(s)\mid \overline{L}_m,\overline{A}_{m-1}=\overline{0},A_m\right]$
does not depend on $A_m$.
\end{assumption}
Assumption~\ref{HH} is a homoscedasticity assumption, because it states that a conditional covariance does not depend on $A_m$. Because of Theorem~\ref{mim} and no unmeasured confounding assumption~\ref{intconf}, homoscedasticity assumption~\ref{HH} is equivalent to
\begin{equation}\label{HHeq}
E\left[H(k)H(s)\mid \overline{L}_m,\overline{A}_{m-1}=\overline{0},A_m\right]=E\left[H(k)H(s)\mid \overline{L}_m,\overline{A}_{m-1}=\overline{0}\right].
\end{equation}
Assumption~(\ref{HHeq}) is not far-fetched: under assumption of no unmeasured confounding~\ref{intconf}, because of Theorem~\ref{mim}, the conditional expectation given $\overline{L}_m,\overline{A}_{m-1}=\overline{0},A_m$ of the two factors $H(k)$ and $H(s)$ does not depend on $A_m$. Assumption~\ref{HH} can be checked empirically by using a preliminary estimator $\tilde{\psi}$ for $\psi_{*}$, regressing the product $H_{\tilde{\psi}}(k)H_{\tilde{\psi}}(s)$ on $\overline{L}_m,\overline{A}_m$, and investigating whether parameter(s) describing the dependence on $A_m$ are equal to $0$.

Rank preservation holds if $H(k)$ does not just mimic $Y_{k}^{(\emptyset)}$ as described in Theorem~\ref{mim}, but $H(k)$ is equal to $Y_{k}^{(\emptyset)}$. Assumption of no unmeasured confounding~\ref{intconf} could be extended, without loss of meaning, to
$\bigl(Y_{k}^{(\emptyset)},Y_{s}^{(\emptyset)}\bigr)\cip A_m\mid \overline{L}_m,$ $\overline{A}_{m-1}=\overline{0}$.
Under this formulation of no unmeasured confounding and rank preservation, equation~(\ref{HHeq}) and thus Assumption~\ref{HH} are immediate. Unfortunately, rank preservation is a very strong assumption, which we do not wish to make. For a discussion see e.g.\ \cite{Enc} or \cite{MimAn}.

The following theorem describes the optimal estimator in an example:
\begin{theorem}\label{qoptHH}\emph{(OPTIMAL ESTIMATOR).}
Suppose that 
\begin{eqnarray*}
\lefteqn{\gamma^{m}_{k,\psi}\left(\overline{L}_m\right)=\left(\psi_1+\psi_2m+\psi_3m^2+\psi_4k\right)(k-m)}\\
&&\;\hspace*{2.5cm}+\left(\psi_5k^2+\psi_6(k-m)+\psi_7lvl_m\right)(k-m)+\psi_8 lvl_m
\end{eqnarray*}
(compare with Example~\ref{exART}).
Suppose that Assumption~\ref{HH} holds. For $k>m$, define
\begin{equation}
\vec{\Delta}_m(k)
=\left(\begin{array}{c}k-m-E\left[Tr(m,k)\mid \overline{A}_{m}=\overline{0},\overline{L}_m\right]\\
m\left(k-m\right)-E\left[T  Tr(m,k)\mid \overline{A}_{m}=\overline{0},\overline{L}_m\right]\\
m^2\left(k-m\right)-E\left[T^2  Tr(m,k)\mid \overline{A}_{m}=\overline{0},\overline{L}_m\right]\\
k\left(k-m-E\left[Tr(m,k)\mid \overline{A}_{m}=\overline{0},\overline{L}_m\right]\right)\\
k^2\left(k-m-E\left[Tr(m,k)\mid \overline{A}_{m}=\overline{0},\overline{L}_m\right]\right)\\
\left(k-m\right)^2-E\left[Tr^2(m,k)\mid \overline{A}_{m}=\overline{0},\overline{L}_m\right]\\
lvl_m(k-m)-E\left[lvl_T  Tr(m,k)\mid \overline{A}_{m}=\overline{0},\overline{L}_m\right]\\
lvl_m-E\left[lvl_T  A_{k-1}\mid \overline{A}_{m}=\overline{0},\overline{L}_m\right]\end{array}\right)\label{Delta}
\end{equation}
with $Tr(m,k)$ the number of treated times between time $m$ and time $k$, and
\begin{equation*}
\vec{cov}_m\left[HI_8\mid \overline{L}_m,\overline{A}_{m-1}=\overline{0}\right]=\left(\begin{array}{ccccc}\Gamma^m_{min,min}I_8&\Gamma^m_{min,min+1}I_8& & &\Gamma^m_{min,max}I_8\\
\Gamma^m_{min+1,min}I_8& & & \\
 & \\
 & \\
\Gamma^m_{max,min}I_8& & & &\Gamma^m_{max,max}I_8\end{array}\right)
\end{equation*}
($I_8$ the $8\times 8$ identity matrix), with
$\Gamma^m_{k,s}={\rm cov}\left[H(k),H(s)\mid \overline{L}_m,\overline{A}_{m-1}=\overline{0}\right]$.
If $\vec{q}^{\;opt}$ satisfies 
\begin{equation*}
\left(\begin{array}{c}\vec{\Delta}_m\left(min\right)\\ \vec{\Delta}_m\left(min+1\right)\\ \\ \\ \vec{\Delta}_m\left(max\right)\end{array}\right)=\vec{cov}_m\left[HI_8\mid \overline{L}_m,\overline{A}_{m-1}=\overline{0}\right]\left(\begin{array}{c}\vec{q}^{\;opt,min}_m\left(\overline{L}_m\right)\\
\vec{q}^{\;opt,min+1}_m\left(\overline{L}_m\right)\\ \\ \\\vec{q}^{\;opt,max}_m\left(\overline{L}_m\right)\end{array}\right),
\end{equation*}
then $P_n\bigl(G^*\bigl(\psi,\hat{\theta},\vec{q}^{\;opt}\bigr)\bigr)=0$ leads to an
optimal estimator $\hat{\psi}$: any other estimator in our class, solving  $P_n\bigl(G^*\bigl(\psi,\hat{\theta},\vec{q}\,\bigr)\bigr)=0$ for some $\vec{q}$, leads to an asymptotic variance that is at least as large as the asymptotic variance of $\hat{\psi}$.
\end{theorem}
The class of estimators solving estimating equations of the form\\ $P_n\bigl(G^*\bigl(\psi,\theta_*,\vec{q}\,\bigr)\bigr)=0$ leads to the smallest asymptotic variance of all estimators
considered in this article (Theorems~\ref{1or3} and~\ref{3or1e}), and
estimating $\theta_{*}$ with pooled logistic regression does not
change the asymptotic variance (Theorem~\ref{thetahat}). $\vec{q}^{\;opt}$
from Theorem~\ref{qoptHH} therefore leads to
the smallest possible asymptotic variance.

\begin{theorem}\label{qoptHH2}\emph{(OPTIMAL ESTIMATOR).}
Extending Theorem~\ref{qoptHH} to different treatment effect models
$\gamma$ can be done as follows. When choosing simpler models for
$\gamma_\psi$, delete the corresponding rows in equation (\ref{Delta}) for
$\vec{\Delta}_m(k)$ and replace the $8$ in $I_8$ by the number of remaining parameters. For more complicated
or different models, notice that the first entry in each row of
$\vec{\Delta}_m(k)$ corresponds to the second entry in each row but with
$A_m=0$ replaced by $A_m=1$. In addition, the model for $\gamma$ of
Theorem~\ref{qoptHH} can easily be generalized to contain similar
terms depending on other covariates; the optimal estimator then
follows similar to Theorem~\ref{qoptHH}.
\end{theorem}

The term $\Gamma^m_{k,s}$ and the double robustness term
$E\left[H_\psi(k)\mid \overline{L}_{m},\overline{A}_{m-1}=\overline{0}\right]$
are fixed functions of $\overline{L}_m$, but they may depend on
$\psi$. We will use an initial estimate $\tilde{\psi}$, doubly robust but not optimal, in place of $\psi_{*}$ in
$\Gamma^m_{k,s}$ and $E\left[H_\psi(k)\mid \overline{L}_{m},\overline{A}_{m-1}=\overline{0}\right]$. If the treatment initiation model $p_\theta$ is correctly specified, the estimating equations are unbiased
for any fixed value of $\tilde{\psi}$. In our application (Section~\ref{results}), a candidate for $\tilde{\psi}$
is motivated by the estimator which is optimal among those with $\vec{q}_m^{\;k}$
only non-zero for $k=m+12$.
This can be shown to lead to $q$ satisfying
\begin{equation*}
\vec{\Delta}_{m}(m+12)={\rm
  Var}\left[H(m+12)\mid \overline{L}_{m},\overline{A}_{m-1}=\overline{0}\right]q^{\;m+12}_{m}\left(\overline{L}_{m}\right).
\end{equation*}
Additionally, replacing the conditional variance of $H(m+12)$ by a working identity
covariance matrix gives $\tilde{q}_{m}^{\;m+12}=\vec{\Delta}_{m}(m+12)$. This
leads to valid estimates $\tilde{\psi}$: stacking the estimating
equations with this $\tilde{q}$ with estimating equations for the
parameters in $\vec{\Delta}_m$, results in unbiased estimating equations
(Theorem~\ref{see}).

The (generalized) inverse of $\vec{cov}_m\left[HI_8\mid \overline{L}_m,\overline{A}_{m-1}=\overline{0}\right]$ equals  the (generalized) inverse of the conditional covariance matrix of $H$ with each entry replaced by the entry times $I_8$.

The optimal estimator requires estimating conditional expectations. In small samples this may be an issue, but as for many efficient estimators (see e.g.~\citeauthor{Newey1993} \citeyear{Newey1993}, \citeyear{Newey1990}, \citeauthor{Hahn} \citeyear{Hahn}, or \citeauthor{Tsiatis} \citeyear{Tsiatis}), it does not lead to a larger asymptotic variance if all models are correctly specified:
\begin{theorem}\label{prelest}
Suppose $\tilde{\psi}_2$ is a preliminary estimator of $\psi_{*}$ which is the result of unbiased estimating equations $P_n\tilde{G}(\psi_2)=0$, $\hat{\theta}$ is an estimator of $\theta_{*}$ from a correctly specified pooled logistic regression model with estimating equations $P_nU(\theta)=0$, and $E_{\xi}\left[H_{\psi}(k)\mid \overline{L}_m,\overline{A}_{m-1}=\overline{0}\right]$ and $\vec{q}_{\psi,\xi}^{\;opt}$ are parameterized by $\xi$, with $\xi_{*}$ the true $\xi$ when $\psi=\psi_{*}$, which can be estimated using estimating equations $P_n J(\xi,\psi)=0$ with $EJ\left(\xi_{*},\psi_{*}\right)=0$. Then, under regularity conditions, solving $\hat{\psi}$ from the unbiased estimating equations
\begin{equation*}P_n\left(\begin{array}{cccc}
G^*(\psi,\theta,\vec{q}^{\; opt}_{m,\psi_2,\xi})\;&\tilde{G}(\psi_2)\;&J(\xi,\psi_2)\;&U(\theta)\end{array}\right)=0,
\end{equation*}
results in the same asymptotic variance for $\hat{\psi}$ as using the
true (but unknown) $\vec{q}^{\;opt}$ from Theorem~\ref{qoptHH} in the
estimating equations\\ $P_n\bigl(G^*\bigl(\psi,\hat{\theta},\vec{q}^{\;opt}\bigr)\bigr)=0$.
\end{theorem}
Solving these estimating equations simultaneously leads to the same estimator $\hat{\psi}$ as plugging in $(\tilde{\psi}_2,\hat{\theta},\hat{\xi})$ into the estimating equations for $\psi_{*}$ and then solving for $\hat{\psi}$. Theorem~\ref{prelest} implies that, if all models are correctly specified, the resulting $\hat{\psi}$ is optimal within the classes studied here.

In practice, instead of estimating $\Gamma^m_{k,s}={\rm
  cov}\left[H(k),H(s)\mid \overline{L}_m,\overline{A}_{m-1}=\overline{0}\right]$,
one may choose to use a so-called working covariance matrix, and
replace $\Gamma^m$ by, for example, the identity matrix. This can be
compared with working correlation matrices in generalized linear models as in
\cite{Zeger1988}. As for generalized linear models, the resulting estimator is not optimal, but
consistency, asymptotic normality and double robustness are not
affected.

\begin{theorem}\label{bootstrap}\emph{(CONSISTENCY OF THE BOOTSTRAP).} Under regularity conditions, the bootstrap for all estimators above is consistent under the conditions already adopted for consistency and asymptotic normality.
\end{theorem}

\section{Applying these methods: estimation steps}\label{steps}

\subsection{Implementation: general remarks}

Section~\ref{steps} details the implementation of the estimators proposed in
this article. The model for $\gamma$ here is
Example~\ref{exART} model~(\ref{modelex}), but the methods can easily be adapted to other treatment effect models. Models that are linear in
$\psi$ are especially attractive because they lead to estimating
equations that are linear in $\psi$ and that are therefore easy to
solve. SAS 9.1.3 (SAS Institute Inc., Cary, North Carolina, USA) was
used for all analyses. The SAS code is available from the author.

As in \cite{Victor}, we adopt a pooled logistic regression model for the treatment prediction model, $p_\theta(m)$. As the outcome $Y_{k}$ one could choose either the CD4 count itself or the CD4 count increase between month $m$ and month $k$. From
Definition~\ref{gamma} of the treatment effect, the same quantity is estimated
whether $Y_{k}$ is the CD4 count or the CD4 count increase: subtracting CD4$_{m}$ from the outcome
CD4$_{k}$ affects both $Y^{(m)}_{k}$ and $Y^{(\emptyset)}_{k}$
the same way, so the CD4$_m$ terms in $\gamma$ cancel. The CD4 count increase
between month $m$ and month $k$ likely reflects more than the
CD4 count at month $k$ the effect of treatment taken between month
$m$ and month $k$, i.e.\ less noise is expected. Thus, the CD4 count increase between month $m$ and month $k$ was used as the outcome for all estimators that are not doubly robust. For doubly
robust estimators, subtracting CD4$_{m}$ from the outcome
CD4$_{k}$ does not affect the point estimates, since the same
covariate is included in $\overline{L}_m$; subtracting
CD4$_{m}$ from CD4$_{k}$ would affect both $Y_{k}$
and $E\left[H(k)\mid \overline{L}_m,\overline{A}_{m-1}=\overline{0}\right]$ in
the same way, so the CD4$_m$ terms in the estimating equations cancel.


\subsection{A preliminary estimator $\tilde{\psi}$}\label{sprel}

As a preliminary estimator $\tilde{\psi}$, necessary to implememt the optimal
estimator, we used a doubly robust version of the estimator from \cite{Victor}.
This choice of $\vec{q}_m^{\;k}$ is the same as in Theorem~\ref{qoptHH2}, but
with $q_m^k=0$ for $k\neq m+12$, and with the
$\vec{cov}_m\left[HI_3\mid \overline{L}_m,\overline{A}_{m-1}=\overline{0}\right]$
replaced by identity matrices. Under a homoscedasticity condition, this
choice of $q$ is optimal within the class of estimating equations with
$\vec{q}_m^{\;k}=0$ for $k\neq m+12$. In the models for $q$,
$E[Tr(m,m+12)\mid \overline{L}_{m},\overline{A}_{m}=\overline{0}]$, etc., we
first estimated $pr(Tr(m,m+12)\neq
0\mid \overline{L}_{m},\overline{A}_{m}=\overline{0})$ using logistic
regression, and then, conditional on $Tr(m,m+12)\neq 0$ and
$\overline{A}_{m}=\overline{0}$, regressed $Tr(m,m+12)$ on the
covariates $\overline{L}_{m}$. In the presence of censoring, we
restricted the regressions to patients still in follow-up at month
$m+12$.
Misspecification of $q$ of the preliminary estimator does not affect
asymptotic optimality or double robustness of the optimal estimator
that makes use of it.  In finite samples it may affect the
variance.

With this treatment effect model, $H(m+12)=Y_{m+12}-(\psi_1+\psi_2T+\psi_3T^2) Tr(m,m+12)$,
and to estimate
$E\left[H(m+12)\mid \overline{L}_{m},\overline{A}_{m-1}=\overline{0}\right]$,
we considered each term in this expression separately, leaving in $\psi$.\\
For $E[Tr(m,m+12)\mid \overline{L}_{m},\overline{A}_{m-1}=\overline{0}]$, we
used the same approach as for
$E[Tr(m,m+12)\mid \overline{L}_{m},\overline{A}_{m}=\overline{0}]$. In the
presence of censoring, we used Inverse Probability of Censoring Weighting (IPCW, see e.g.~\citeauthor{RRZ} \citeyear{RRZ}) starting at month $m$.
With $C_p=0$ indicating a patient was uncensored at month $p$, these weights are
\begin{equation*}
W_{m,k}=\left\{\prod_{p=m+1}^{k}pr\left(C_p=0\mid \overline{L}_{p-1},\overline{A}_{p-1},\overline{C}_{p-1}=\overline{0}\right)\right\}^{-1}.
\end{equation*}

This procedure leads to estimating equations that are linear in $\psi$ and
thus easy to solve.

\subsection{The optimal estimator}\label{sopt}

The mimicking outcome $H(k)$ was first estimated by $H_{\tilde{\psi}}(k)$. Theorem~\ref{mim} facilitates specification of a model for $E[H(k)\mid \overline{L}_m,\overline{A}_{m-1}=\overline{0}]$. Linear
regression was used to estimate $E[H(k)\mid \overline{L}_m,\overline{A}_{m-1}=\overline{0}]$. In the
presence of censoring, we used IPCW, as for the preliminary estimator. This leads to an estimate
$\hat{E}[H(k)\mid \overline{L}_m,\overline{A}_{m-1}=\overline{0}]$ based on
estimating equations. Optimality depends on
correct specification of $E[H(k)\mid \overline{L}_m,\overline{A}_{m-1}=\overline{0}]$, but
if the treatment initiation model $p_\theta$ is correctly specified, consistency and asymptotic normality do not (because of double robustness). In the model for
$\Delta$, $E[Tr(m,k)\mid \overline{L}_m,\overline{A}_{m}=\overline{0}]$
etcetera, the same approach as in Section~\ref{sprel} was used:
first $pr(Tr(m,k)\neq
0\mid \overline{L}_m,\overline{A}_{m}=\overline{0})$ was estimated using logistic
regression; then, conditional on $Tr(m,k)\neq 0$ and
$\overline{A}_m=\overline{0}$, $Tr(m,k)$ was regressed on the
covariates $\overline{L}_m$. In the presence of censoring, the
regression was restricted to patients still in follow-up at month
$k$. In the simulations, $\vec{cov}_m\left[HI_3\mid \overline{L}_m,\overline{A}_{m-1}=\overline{0}\right]$ (Theorem~\ref{qoptHH2}) does not depend on $\overline{L}_m$. In the
application, a working model not depending on $\overline{L}_m$
was used for $\vec{cov}_m\left[HI_3\mid \overline{L}_m,\overline{A}_{m-1}=\overline{0}\right]$,  similar to a working covariance matrix. $\Gamma_{k,s}^m$ was estimated by the empirical
average over all patients of\\
$\left\{H_{\tilde{\psi}}(k)-\hat{E}[H(k)\mid \overline{L}_m,\overline{A}_{m-1}=0]\right\}\left\{H_{\tilde{\psi}}(s)-\hat{E}[H(s)\mid \overline{L}_m,\overline{A}_{m-1}=0]\right\}$.\\ Alternatively, also after plugging in $\tilde{\psi}$, we
could have used techniques similar to GEE to estimate a
working covariance matrix. In the presence of censoring, we added
$\overline{C}_{{\rm max}(k,s)}=\overline{0}$ to the conditioning
event.
Misspecification of $q$ (through $\Delta$ or $\vec{cov}_m\left[HI_3\mid \overline{L}_m,\overline{A}_{m-1}=\overline{0}\right]$) leads to a suboptimal
estimator, but (Theorem~\ref{prelest}) does not affect double robustness, because all specifications leading to $q$ only depend
on $\overline{L}_m$ and parameters solving estimating equations.

This procedure leads to estimating equations that are linear in $\psi$ and
thus easy to solve.

\subsection{For comparison, a naive choice}\label{snaive}

For comparison, we implemented two non-doubly-robust estimators, based
on Theorem~\ref{see} and not using the optimality theory
developed here. For these Theorem~2-based estimators, since in our application (Section~\ref{results}) interest lies in the effect of one year of treatment, $\vec{q}_m^{\;k}=0$
for $k\neq m+12$, and the $\vec{q}_{m}^{\;m+12}$ were as follows, with CD4$_m$ the CD4 count at month $m$, injdrug an indicator of whether the patient ever injected drugs at or before the first visit, lvl$_m$ the log$_{10}$ viral load at month $m$, and firstvisit$_m$ an indicator for whether month $m$ was the month of the first visit. In the
simulations, we chose
\begin{equation}\label{naive1sim}
\;\hspace{2cm}\vec{q}_{m}^{\;m+12}=\bigl(\begin{array}{ccc}{\rm CD4}_m\; & m\; & {\rm injdrug}\end{array}\bigr)^{\top}\vspace*{-0.1cm}
\end{equation}
\begin{equation}\label{naive2sim}
\text{and}\hspace{2cm} \vec{q}_{m}^{\;m+12}=\bigl(\begin{array}{ccc}{\rm CD4}_m\;& {\rm injdrug}\; & {\rm CD4}_6\end{array}\bigr)^\top.
\end{equation}
In the data application, we chose
\begin{equation}\label{naive1}
\;\hspace{2cm}\vec{q}_{m}^{\;m+12}=\bigl(\begin{array}{ccc}\sqrt{{\rm CD4}}_m \;& m\; & {\rm lvl}_m \end{array}\bigr)^{\top}
\end{equation}
\begin{equation}\label{naive2}
\text{and}\hspace{2cm} \vec{q}_{m}^{\;m+12}=\bigl(\begin{array}{ccc}\sqrt{{\rm CD4}}_m\; & {\rm lvl}_m\; & {\rm firstvisit}_m\end{array}\bigr)^{\top}.
\end{equation}

\section{Simulations}\label{sim}

We simulated data with monthly visits, and based choices for the simulated
data on the AIEDRP data on HIV infected patients, described in Section~\ref{results}. We used an auto-regressive
model for the course of the CD4 count, which may be more realistic in
months 6-30 than before month $6$, given the different behavior of CD4 counts in
the first 6 months since HIV infection (Web-Appendix~C).
Therefore, we simulated data in months 6-30, and estimated the effect of
treatment initiation in months 6-18. Simulations are
detailed in Web-Appendix~C.


We simulated two scenarios: 1.: 1000 datasets
with 1000 observations each, and 2.: 500 datasets with 5000
observations each. We fitted model (\ref{modelex}) with 2 parameters, $(\psi_1,\psi_2)$, and 3 parameters, $(\psi_1,\psi_2,\psi_3)$. In these simulations, the true $\psi_3$ equals
$0$. Since in our application interest focused on the effect of one year of treatment, we restricted the estimating equations to $k=m+1,\ldots,m+12$ so as to rely less on model specification. Table~1 shows the root mean squared errors of
the different estimators described in Sections~\ref{opt} and~\ref{steps}. Web-Appendix~C.2 describes the models fitted for the nuisance parameters.

Table~1 shows great improvements from applying our theory in
comparison with a naive choice of estimating equations (Table~1 estimators
1a.\ and 1b., described in Section~\ref{snaive}). Choosing $q$ without using optimality theory can lead to useless inference. Making the estimator doubly robust (estimator 3.) results in improvements of the mean squared errors in the 3-parameter model. Not restricting the analysis to $k=m+12$ (estimators 4.\ and 5.) results in improvements overall, and our estimator with optimal asymptotic properties (estimator 5.) performs best in the simulations.

We also calculated the $2.5\%$ and $97.5\%$ quantiles (over the datasets) of the estimated parameters for estimators 2., 3., 4., and 5. For all parameters, the truth was in between these quantiles. 

For $n=1000$, we also investigated choosing sparser models (Web-Appen-dix~C.2) in the expressions for the prediction of treatment duration (for $q$ and for the doubly robust term). Choosing sparser models made little difference for estimators 3., 4., and 5 in Table~1. For estimator~2.\ in Table~1, sparser models for $q$ led to substantially larger mean squared errors (results not shown).

Figure~1 shows the results for the datasets with 1000
observations. Figure~1 compares the performance of estimators of the
effect of treatment on the 1-year increase in the CD4s count due to
treatment initiated at the different months. For example, for month 11,
this is the root mean squared error for the expected difference in CD4 count at month
23=12+11 between 1.\ initiating ART at month 11 versus 2.\ never
initiating ART. Figure~1 does not include the naive estimators (1a.\ and
1b.\ in Table~1), since those perform much worse and incorporating them
makes a comparison of the other estimators impossible. The same
pattern appears as in Table~1, with the optimal estimator performing best. 


\section{The effect of ART in HIV infected patients during acute and early HIV infection}
\label{results}

ART is the standard of care for HIV infection. Guidelines regarding ART initiation have been changing,
with patients initiating ART earlier (\citeauthor{Thom} \citeyear{Thom}, \citeauthor{panel5} \citeyear{panel5}, \citeauthor{WHO2016} \citeyear{WHO2016}), often even at the time of diagnosis, especially in the developed world. 


The effect of ART in the early and acute stages of
HIV infection was studied in \cite{Victor} using a preliminary version of
our article; we aim to estimate the effect of ART with more
robustness and precision. Given the limited data on early HIV infection, this is a timely question in HIV research. We estimate how the effect on immune reconstitution of initiating one year of
ART depends on the time between the estimated date of HIV infection and
ART initiation. The effect on immune
reconstitution of one year of ART initiated $m$ months after
infection, is measured as the CD4 count at month $m+12$ with ART
initiated at month $m$ versus the CD4 count at month $m+12$ without ART. Of particular interest is the effect of one year
of treatment initiated at month $m$, given past covariate history
$\overline{l}_m$; that is, 
$\gamma^{m}_{m+12}\left(\overline{l}_m\right)$.

The results of our investigation are important, since ART initiation
soon after infection may not only improve a patient's own outcomes, but
also reduces the risk that a patient's HIV infection is spread to
others (\citeauthor{Cohen} \citeyear{Cohen}, \citeauthor{Granich} \citeyear{Granich}, \citeauthor{Contr} \citeyear{Contr}). Furthermore, the
CDC is encouraging HIV testing (\citeauthor{CDC} \citeyear{CDC}), leading to earlier HIV
diagnoses, so more treatment initiation decisions need to be made
during early and acute infection. There is not a lot of evidence of clinical benefit for initiating ART this early, with likely a relatively small number of patients in the START trial (\citeauthor{StartTrial} \citeyear{StartTrial}) in acute or early infection at baseline. 
Our investigations shed light on the effect of efforts to diagnose HIV early, if early diagnosis is combined with immediate
ART initiation as currently recommended.

We apply our estimators to the observational AIEDRP
(Acute Infection and Early Disease Research Program) Core01 data, using
data on 1762 HIV infected patients diagnosed during acute and early
infection (\citeauthor{Hecht} \citeyear{Hecht}). Dates of infection are estimated using an
algorithm that incorporates clinical and laboratory data (\citeauthor{Hecht} \citeyear{Hecht},
\citeauthor{Smith} \citeyear{Smith}). \cite{Victor} showed that in the AIEDRP, ART use depends
on covariates such as the current CD4 count that are prognostic for the
outcome CD4 count; and that this leads to substantial confounding by indication.

In this HIV application,
$K+1$ is 24 months. $0$ is the estimated date of
HIV infection, although visits may be missed during follow-up and
particularly in the earliest months of infection. To account for missed visits, we include the visit
pattern (in which months a visit took place) as a measured
covariate. $Y$, $A$, and $L$ are measured at multiple time points
that vary across patients. For $L$ we use the average measurement
within the given month; if this is missing at month $m$, $L_m$ is
coded as missing, a possible covariate value. $A_m$ cannot be
missing because we assume that treatment can only start at visits and
then it is always recorded. We impute missing data on the outcomes $Y_k$,
after the first visit and until censoring, by interpolation, except
for visits just prior to onset of treatment, for which we carry the
last observation forward, in order to avoid using post-treatment information to impute outcomes prior to treatment. We assumed:
\begin{assumption} \label{modelappl} \emph{(Parameterization of coarse structural nested mean model).} For $k=(m+1)\vee 12,\ldots,(m+12)\wedge (K+1)$, suppose that 
\begin{equation*}
\gamma^{m}_{k,\psi}\left(\overline{l}_m\right)=
\left(\psi_{*1}+\psi_{*2}m+\psi_{*3}m^2\right)(k-m)1_{\left\{k>m\right\}},
\end{equation*} with $(k-m)$ the treatment duration from
month $m$ to month $k$.
\end{assumption}
Since interest in this HIV application lies in the effect of 1 year of treatment, $\gamma^{m}_{m+12}\left(\overline{l}_m\right)$, Assumption~\ref{modelappl} suffices. Specifying $\gamma_k^m$ for $k<12$ and $k>m+12$ might lead to greater precision, but increases the risk of model misspecification. The restriction to these values of $k$ implies that in assumption of no unmeasured confounding~\ref{nuc}, $k$ can be similarly restricted. For estimation, every sum over $k$ is then also restricted to $k=(m+1)\vee 12,\ldots, (m+12) \wedge (K+1)$. 

Because treatment is assumed to only change at visit times, the estimating equations include $1_{visit}(m)$, an indicator of whether a visit took place at time $m$. For loss to follow-up, we assumed missing at random (\citeauthor{RubinMAR} \citeyear{RubinMAR}) and applied inverse probablity of censoring weighting (\citeauthor{RRZ} \citeyear{RRZ}). Web-Appendix~D describes the nuisance parameter models.

Table~2 provides estimates and bootstrap $95\%$ confidence intervals based on the AIEDRP data, for the same estimators as described in Section~\ref{sim}. Table~2 shows that importantly, all estimators based on the optimality theory
of the current article lead to much narrower $95\%$ confidence intervals than the
estimators based on Theorem~\ref{see} with $q$ chosen in a naive
way as described in Section~\ref{snaive}. Both naive estimators
lead to irrelevant estimators in the 3-parameter model, due to
extremely wide confidence intervals. The same is true for the second
naive estimator in the 2-parameter model. Double robustness does not lead to narrower confidence
intervals in the AIEDRP data (compare estimators 2.\ and 3.). Estimator 5., which would be
optimal without censoring and under homoscedasticity assumption~\ref{HH}, leads to much wider confidence
intervals than Estimator 4. For the AIEDRP data, Estimator 4., which is
similar to the optimal estimator but with working identity covariance
matrices, leads to the narrowest confidence intervals.

Figure~2 compares the performance of estimators of the effect of
ART treatment on the one-year increase in the CD4 count due to ART initiated
at the different months since the estimated date of HIV infection. For
example, for month 11, the quantity in Figure~2 is the estimated expected difference in the CD4 count at month 23=11+12, comparing initiating ART at month 11 versus never initiating ART. To facilitate
the comparison of the other estimators, Figure~2 does not include the
naive estimators of Section~\ref{snaive}, which performed much
worse. As in Table~2, Estimator~4., with a working identity covariance matrix, leads to the best precision.

Table~2 also describes the results of a sensitivity analysis (details in Web-Appendix~D). In this sensitivity analysis, treatment initiation and dropout are modeled using model selection techniques. While model selection in principle invalidates the confidence intervals, this sensitivity analysis indicates that the results are somewhat sensitive to model specification.

The estimated effect of one year of ART initiated during acute and early HIV infection is substantial and significant. It decreases somewhat when the time between the estimated date of infection and ART initiation increases, but this
trend is insignificant.

\section{Discussion}
\label{Disc}

Causal inference methods to analyze observational data require untestable assumptions, so it
is important to analyze longitudinal observational data using different methods and compare the results. This
requires a variety of methods, and the further development of structural nested models, as an
addition to the wider applied class of weighting methods.

Both in the simulation study and in the HIV application, most of the naive coarse structural
nested mean model (coarse SNMM) estimators led to useless inference. This could realistically happen in practice
if the nuisance function $\vec{q}$ is chosen without knowledge of optimality theory, and motivated this
article. In the HIV application, the naive estimators were so imprecise that it was
important to not use a naive preliminary estimator in the first step for the optimal estimator.

Our theory resulted in substantially improved precision of coarse SNMMs. In the simulation study, the
optimal doubly robust estimator resulted in useful inference, and performed best. In the HIV application,
our methods also substantially improved precision. In the HIV application, the best
performance was by a doubly robust estimator related to the optimal estimator, but using
working identity covariance matrices, similar to the identity working covariance matrices approach
in Generalized Estimating Equations (\citeauthor{Tsiatis} \citeyear{Tsiatis}, Section~4.6). The suboptimal behavior of the optimal estimator may
have several causes. It could be due to the combination of limited sample size and censoring.
It could also be due to model misspecification of the nuisance parameter models, especially of
$\vec{cov}_m[HI_3|\overline{L}_m,\overline{A}_{m-1}=0]$. We focused on coarse SNMMs with a time-varying outcome; Web-appendix~B shows how the calculations simplify considerably with an outcome measured at the end of the study.

Inverse probability of treatment weighting of marginal structural models (\citeauthor{MSM1} \citeyear{MSM1}, \citeauthor{MSM2} \citeyear{MSM2}) can only be used to estimate how a treatment effect depends on baseline covariates. It would be interesting to know whether coarse SNMMs provide more precise estimators than marginal structural models when interest lies in this scenario.
Efficiency gains of coarse SNMMs could be expected, because SNMMs use all observed data, whereas if a saturated outcome model is used, marginal structural
models use only the data that is consistent with the specific treatment regime. Our investigation
is the first step in the comparison between coarse SNMMs and marginal structural models: we optimized of the performance of coarse SNMMs.

A promising area for future research are estimation methods when the number of potential
confounders is large. Regularization methods such as LASSO are then the obvious candidates to
estimate the propensity scores and the conditional expectation of H, assuming e.g. that only a
limited number of confounders truly contribute. Since our proposed estimators are doubly robust,
the estimating equations satisfy the orthogonality or immunization condition of \cite{Chernozhukov2015}. This can be seen since, as a function of one of the two nuisance parameter models,
say $\xi$, the estimating function has expectation $0$ for all $\xi$ when the other parameter, say $\theta$, is
held at its true value. Thus, provided that a high-quality regularization method is used, it can be
expected that if both nuisance parameter models are correctly specified, asymptotically correct
inference can be obtained by using the methods proposed in \cite{Chernozhukov2015}. Correct
specification of nuisance parameter models will be more likely for reasonable sample sizes if
regularization is used, allowing for more elaborate candidate models. In addition, Yang and Lok are investigating semiparametric efficiency in the presence of censoring; this will involve investigating semiparametric efficient estimation as in e.g. \cite{Tsiatis} and \cite{Hahn}, but adapted to treatment effects that depend on time-dependent pre-treatment covariates. Targeted maximum likelihood estimation could be another option to investigate semiparametric efficiency.

This article has limitations. As in traditional structural nested mean models, we have assumed
that all confounders are measured. This assumption cannot be tested using the available data.
Subject matter experts have to judge whether the data analyst included all covariates that are
predictive of both treatment initiation and outcome. Moreover, in our HIV application we have
assumed that a simple treatment effect model is correctly specified. \cite{Shu} tested
this assumption, and concluded that the AIEDRP data do not provide evidence that this assumption
is violated; however, this may be due to a limited sample size. Investigating the properties
of the proposed estimator when the treatment effect model is misspecified is an interesting topic
for future research. In addition, our calculations for optimality are restricted to linear treatment
effect models. Future research may involve nonlinear treatment effect models; for this \cite{Newey1990}, although mainly aimed at estimating the effect of point exposures, may be useful. Moreover,
our estimators are likely not optimal in the presence of censoring. They do remain doubly
robust if censoring is missing at random (\citeauthor{RubinMAR} \citeyear{RubinMAR}) and the model for censoring is correctly
specified. In our HIV application, which includes censoring, our methods still led to remarkable
variance reductions. In addition, one or more of the models for the nuisance parameters may
be misspecified. Only the treatment initiation model and the outcome regression model affect consistency and asymptotic normality, and double robustness implies that misspecification of one of these two models preserves consistency and asymptotic
normality. 


We conclude that the precision of estimators for coarse structural nested mean models depends
substantially on the estimating equations chosen. The substantial improvement we found
by choosing optimal estimating equations suggests that the use of optimal estimators may encourage more widespread use of coarse structural nested mean models.

\section*{Acknowledgements}

This work was supported by
the Milton Fund, the Career Incubator Fund from the Harvard
School of Public Health, NSF DMS 1854934, and the National Institutes of Health [grant numbers NIAID R01 AI100762,
R37 51164
AI43638, AI074621, AI106039
and AI036214].
The content is solely the responsibility of the authors and does not necessarily represent the official views of the National Institutes of Health or the National Science Foundation.

I am grateful to the patients who volunteered for AIEDRP, to the\\
AIEDRP study team, and to Susan Little, Davey Smith, and Christy
Anderson for their help and advice in interpreting the AIEDRP
database. I would like to thank Ray Griner for extensive help with the programming in SAS, and Shu Yang for programming the estimators additionally in R.
I would like to thank James Robins and Victor DeGruttola
for insightful and fruitful discussions.

\section*{Supplementary material}

The appendix includes proofs, outcomes measured at the end of the study, further details of the simulation study, and a description of the nuisance parameter models used in the HIV application.

\addcontentsline{toc}{chapter}{Bibliography}
\bibliographystyle{chicago} \bibliography{ref}

\begin{figure}
\centerline{\includegraphics[scale=0.95, angle=0]{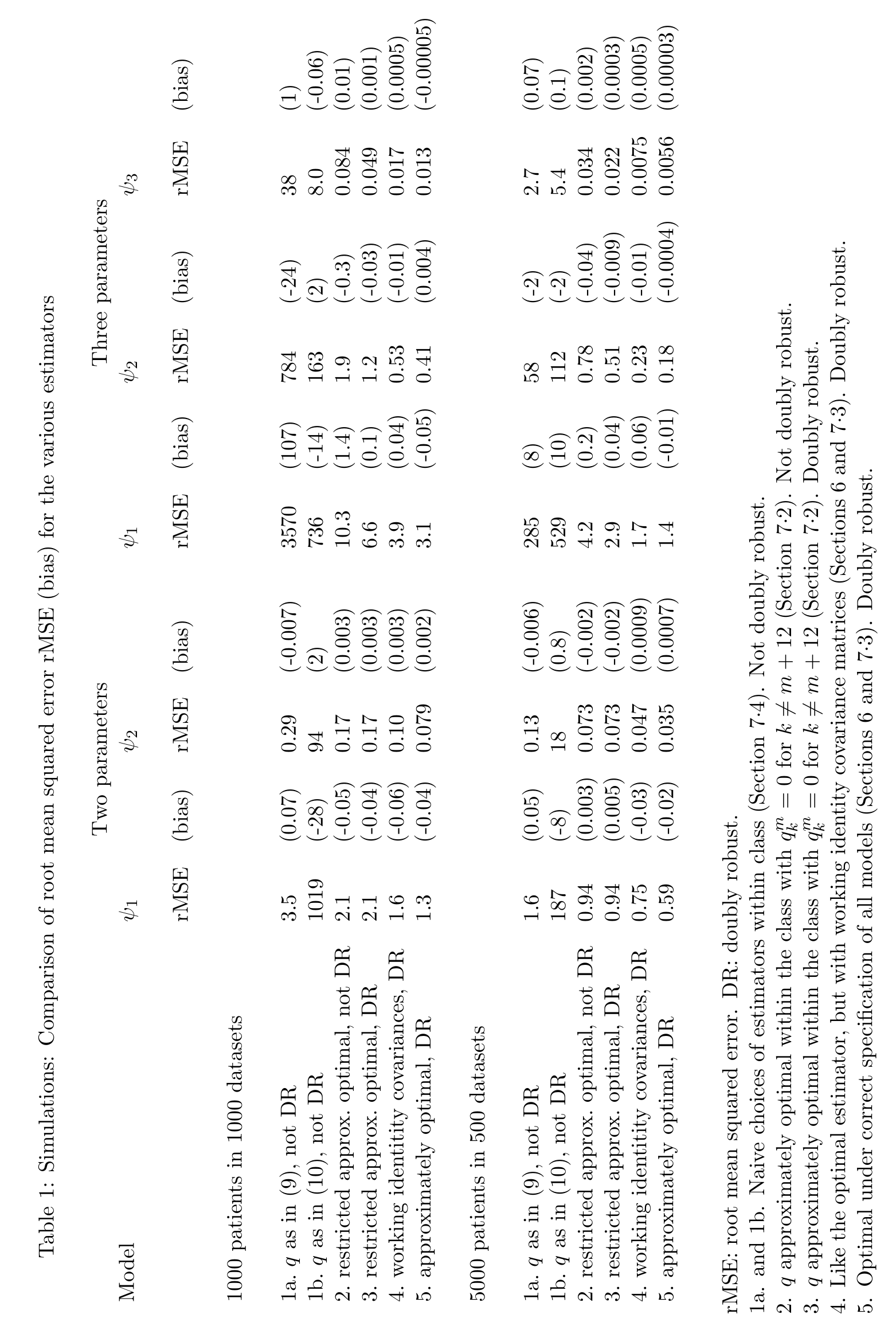}}
\end{figure}

\begin{figure}
\centerline{\includegraphics[scale=1, angle=0]{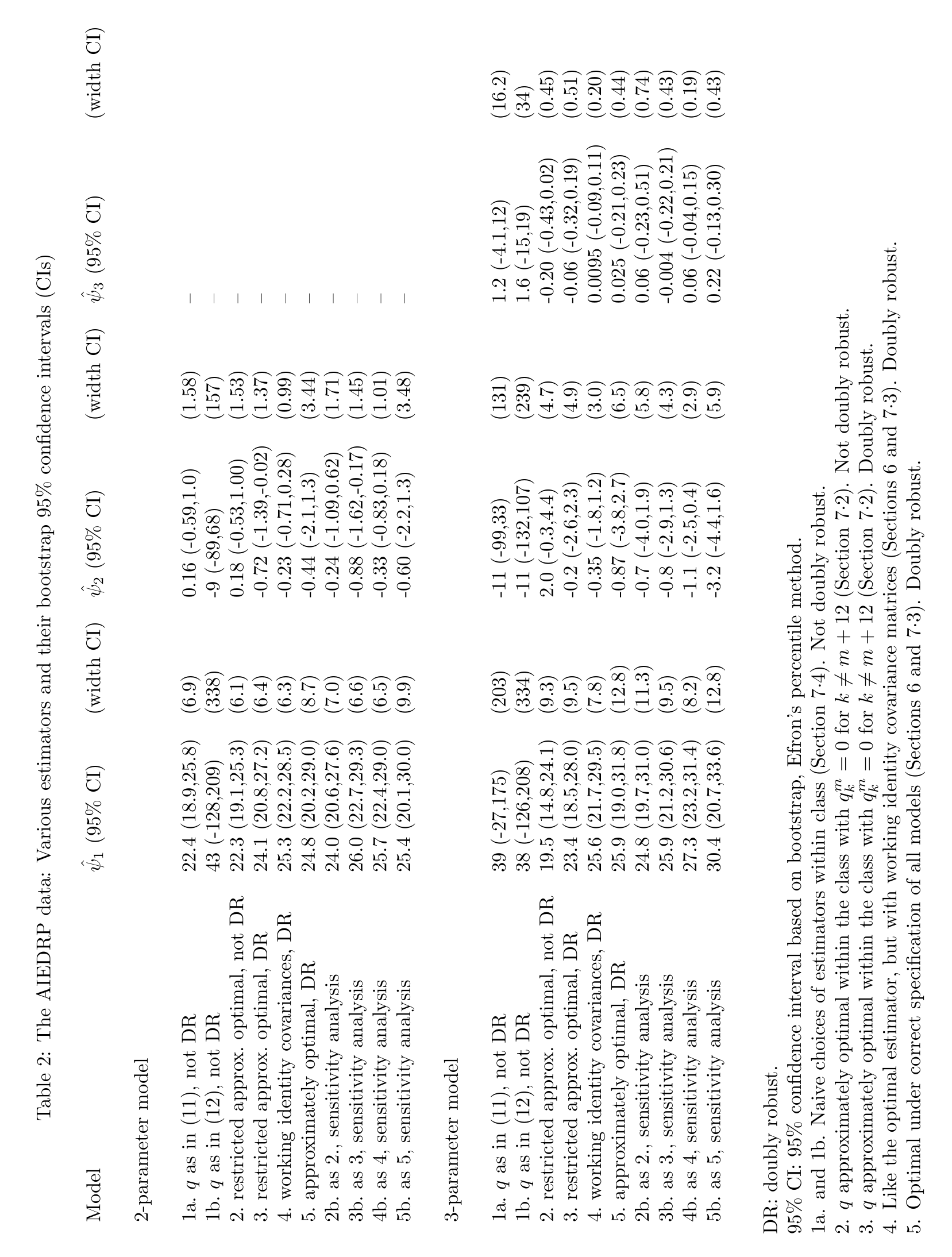}}
\end{figure}

\begin{figure}
\centerline{\includegraphics[scale=1, angle=0]{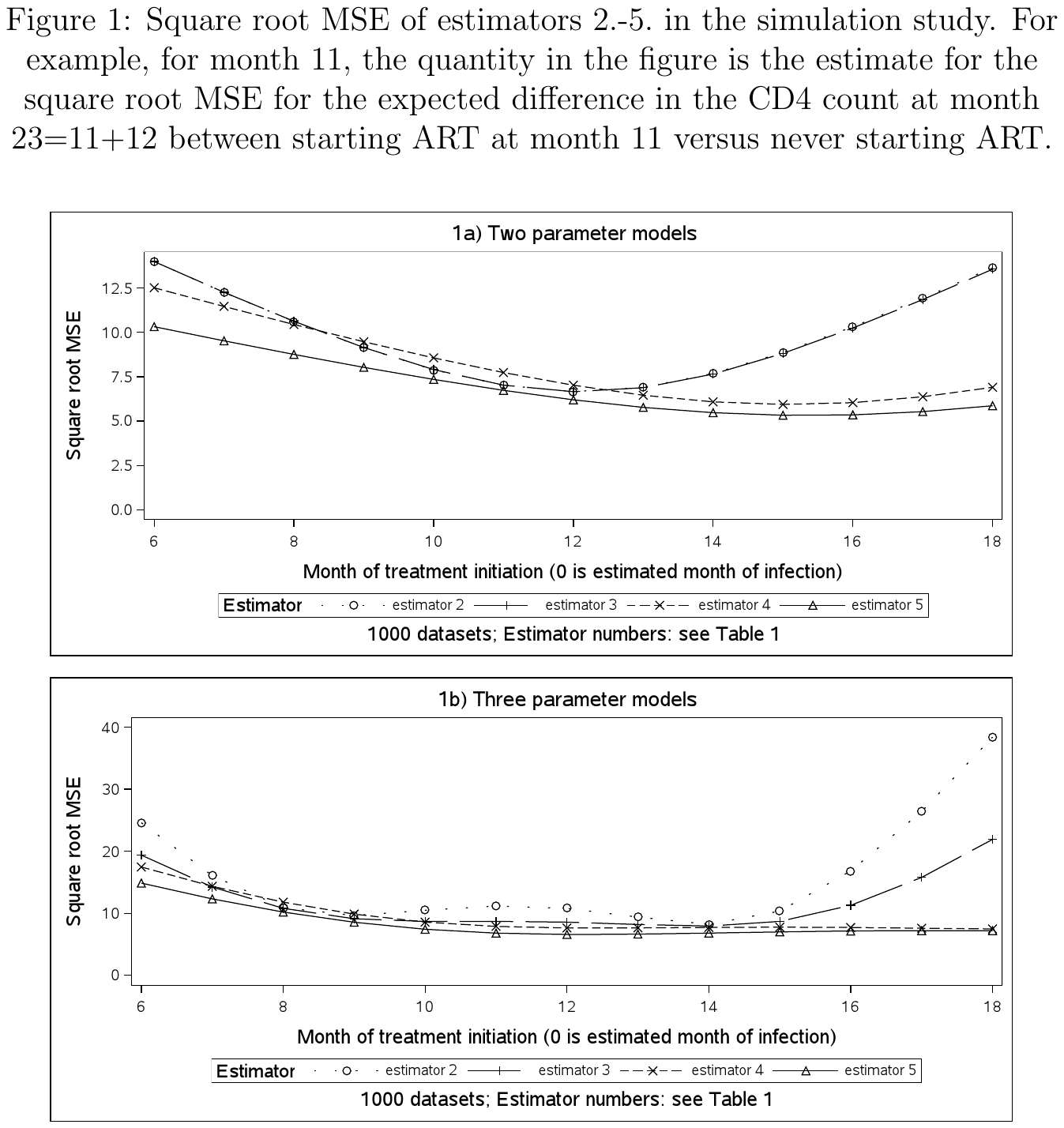}}
\end{figure}

\begin{figure}
\centerline{\includegraphics[scale=1, angle=0]{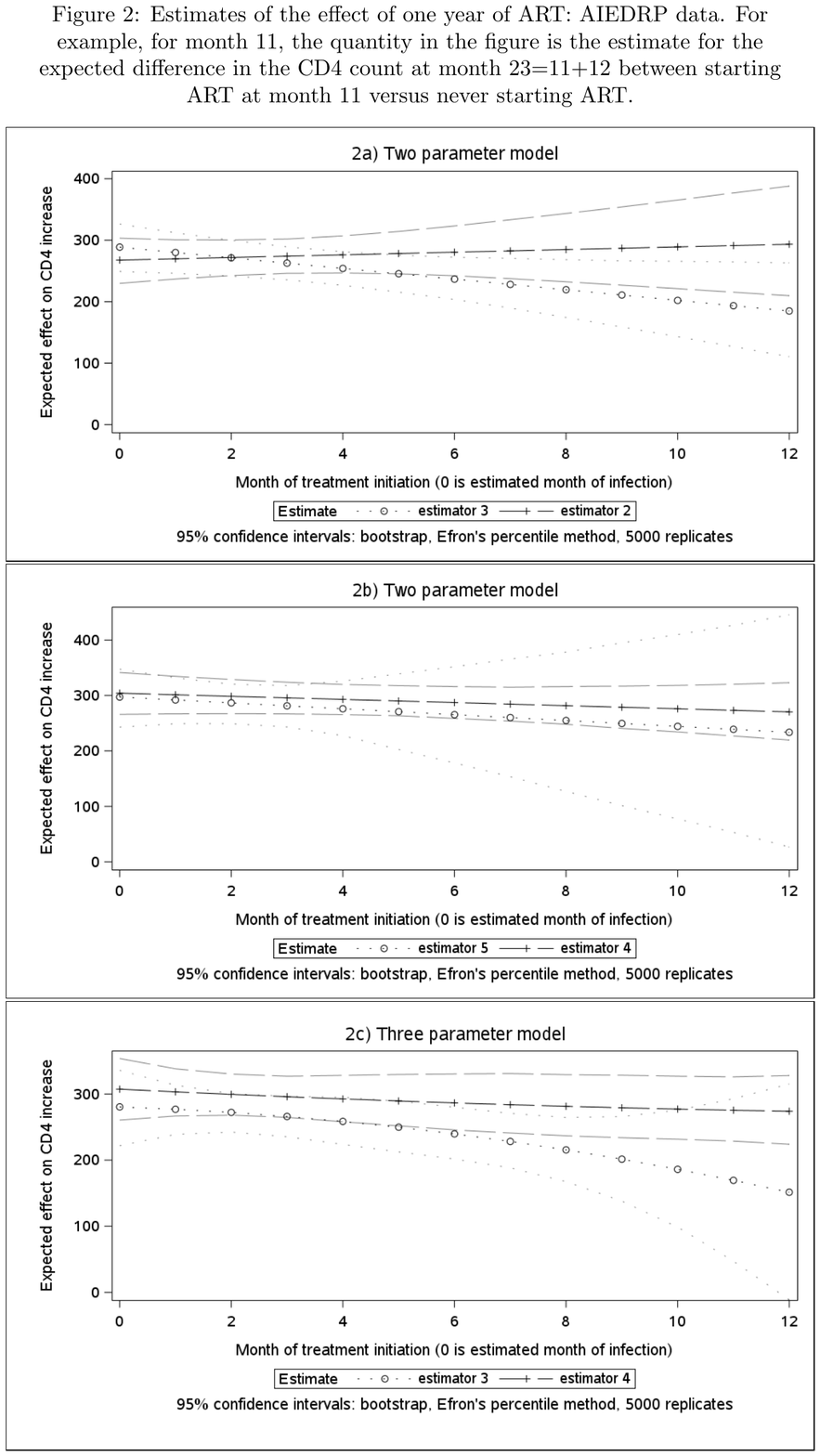}}
\end{figure}

\end{document}